\documentclass[graybox]{svmult}

\usepackage{type1cm}        
\usepackage{makeidx}         
\usepackage{graphicx}        
\usepackage{multicol}        
\usepackage[bottom]{footmisc}

\usepackage{newtxtext}       %
\usepackage[varvw]{newtxmath}       
\usepackage{subcaption}
\usepackage{multirow}

\DeclareMathAlphabet{\pazocal}{OMS}{zplm}{m}{n}

\RequirePackage[ruled,vlined,linesnumbered]{algorithm2e}

\SetCommentSty{xCommentSty}

\definecolor{kothari_Gray}{gray}{0.9}

\newcolumntype{"}{@{\hskip\tabcolsep\vrule width 1pt\hskip\tabcolsep}}

\makeatletter
\newcommand{\thickhline}{\noalign {\ifnum 0=`}\fi \hrule height 1pt
    \futurelet \reserved@a \@xhline
}
\makeatother

\makeindex             

\begin{document}

\title*{A Multigrid Preconditioner for Jacobian-free Newton-Krylov Methods}

\author{Hardik Kothari, Alena Kopani\v{c}\'{a}kov\'{a} and Rolf Krause}
\institute{Hardik Kothari  \at Euler Institute, Universit\`{a} della Svizzera italiana, Switzerland, \email{hardik.kothari@usi.ch}
\and Alena Kopani\v{c}\'{a}kov\'{a} \at  Euler Institute, Universit\`{a} della Svizzera italiana, Switzerland, \email{alena.kopanicakova@usi.ch}
\and Rolf Krause \at  Euler Institute, Universit\`{a} della Svizzera italiana, Switzerland, \email{rolf.krause@usi.ch}
}
\maketitle

\abstract*{In this work, we propose a multigrid preconditioner for Jacobian-free Newton-Krylov (JFNK) methods.
Our multigrid method does not require knowledge of the Jacobian at any level of the multigrid hierarchy.
As it is common in standard multigrid methods, the proposed method also relies on three building blocks: transfer operators, smoothers, and a coarse level solver.
In addition to the restriction and prolongation operator, we also use a projection operator to transfer the current Newton iterate to a coarser level.
The three-level Chebyshev semi-iterative method is employed as a smoother, as it has good smoothing properties and does not require the representation of the Jacobian matrix.
We replace the direct solver on the coarsest-level with a matrix-free Krylov subspace method, thus giving rise to a truly Jacobian-free multigrid preconditioner.
We will discuss all building blocks of our multigrid preconditioner in detail and demonstrate the robustness and the efficiency of the proposed method using several numerical examples.
}

\abstract{In this work, we propose a multigrid preconditioner for Jacobian-free Newton-Krylov (JFNK) methods.
Our multigrid method does not require knowledge of the Jacobian at any level of the multigrid hierarchy.
As it is common in standard multigrid methods, the proposed method also relies on three building blocks: transfer operators, smoothers, and a coarse level solver.
In addition to the restriction and prolongation operator, we also use a projection operator to transfer the current Newton iterate to a coarser level.
The three-level Chebyshev semi-iterative method is employed as a smoother, as it has good smoothing properties and does not require the representation of the Jacobian matrix.
We replace the direct solver on the coarsest-level with a matrix-free Krylov subspace method, thus giving rise to a truly Jacobian-free multigrid preconditioner.
We will discuss all building blocks of our multigrid preconditioner in detail and demonstrate the robustness and the efficiency of the proposed method using several numerical examples.
}

\section{Introduction}
\label{sec:kothari_h_mini_17_intro}
The numerical solution of partial differential equations (PDEs) is often carried out using discretization techniques, such as the finite element method (FEM), and typically requires the solution of a nonlinear system of equations.
These nonlinear systems are often solved using some variant of the Newton method, which utilizes a sequence of iterates generated by solving a linear system of equations.
However, for problems such as inverse problems, optimal control problems, or higher-order coupled PDEs, it can be computationally expensive, or even impossible to assemble a Jacobian matrix.

The Jacobian-free Newton Krylov (JFNK) methods exploit the finite difference method to evaluate the action of a Jacobian on a vector, without requiring the knowledge of the analytical form of the Jacobian and still retain local quadratic convergence of the Newton method.
Even though JFNK methods are quite effective, the convergence properties of the Krylov subspace methods deteriorate with increasing problem size.
Hence, it is desirable to reduce the overall computational cost by accelerating the convergence of the Krylov methods.
To this end, many preconditioning strategies have been proposed in the literature, see e.g.,~\cite{knoll_jacobian-free_2004}.
We aim to employ multigrid (MG) as a preconditioner to accelerate the convergence of the Krylov subspace methods.
Unfortunately, it is not straightforward to incorporate the MG method into the JFNK framework, as the standard implementations of the MG method require either a matrix representation of the Jacobian or an analytical form of the Jacobian.

In this work, we propose a matrix-free geometric multigrid preconditioner for the Krylov methods used within the JFNK framework.
The proposed method exploits the finite difference technique to evaluate the action of Jacobian on a vector on all levels of multilevel hierarchy and does not require explicit knowledge of the Jacobian.
Additionally, we employ polynomial smoothers which can be naturally extended to a matrix-free framework.
Compared to other matrix-free MG preconditioners proposed in the literature, e.g.,~\cite{bastian_matrix-free_2019,davydov_matrix-free_2020,mavriplis_assessment_2002, may_scalable_2015}, our method does not require the knowledge of the analytical form of the Jacobian, and no additional modifications are required in the assembly routine to compute the action of a Jacobian on a vector.

\vspace{0.3cm}

\noindent
\textbf{Jacobian-free Newton-Krylov methods:}
The Newton method is the most frequently used iterative scheme for solving nonlinear problems.
Newton method is designed to find a root $\boldsymbol{x}^{\ast} \in \mathbb{R}^n$ of some nonlinear equation $F(\boldsymbol{x}^{\ast})~=~0$.
The iteration process has the following form:
\[
  \boldsymbol{x}^{(k+1)} = \boldsymbol{x}^{(k)} +\alpha \delta \boldsymbol{x}^{(k)}, \quad \text{for} \  k=0,1,2,\dots,
\]
where~$\alpha > 0$ denotes a line-search parameter and $\delta \boldsymbol{x}^{(k)}$ denotes a Newton direction.
The correction $\delta \boldsymbol{x}^{(k)}$ is obtained by solving the following linear system of equations: ${J(\boldsymbol{x}^{(k)}) \delta \boldsymbol{x}^{(k)} = -F(\boldsymbol{x}^{(k)})}$, where $J(\boldsymbol{x}^{(k)}) = \nabla F(\boldsymbol{x}^{(k)})$.
In the context of this work, we assume that the $F$ is obtained as a gradient of some energy functional $\Psi$, i.e.,~$F(\boldsymbol{x}^{(k)}) \equiv \nabla \Psi(\boldsymbol{x}^{(k)})$.
In this way, the Jacobian $J$ will be a symmetric matrix, which in turn allows us to use a multigrid preconditioner.
In the JFNK methods~\cite{knoll_jacobian-free_2004}, the solution process is performed without explicit knowledge of the Jacobian~$J$.
Instead, the application of a Jacobian to a vector is approximated using the finite difference scheme, given as
{\(
  J(\boldsymbol{x}^{(k)}) \boldsymbol{u} \approx \frac{F(\boldsymbol{x}^{(k)} +\epsilon \boldsymbol{u}) - F(\boldsymbol{x}^{(k)})}{\epsilon},
\)}
where we choose $\epsilon = \dfrac{1}{n\|\boldsymbol{u}\|_2} \sum_{i=1}^n \sqrt{\varepsilon_p} (1+|x^{(k)}_i|)$ and~$\varepsilon_p$ denotes the machine precision.
The value of the finite difference interval $\epsilon$ is chosen, such that the approximation of the Jacobian is sufficiently accurate and is not spoiled by the roundoff errors.

\section{Matrix-free Multigrid Preconditioner}
\label{sec:kothari_h_mini_17_JFMG}
The multigrid method is one of the most efficient techniques for solving linear systems of equations stemming from the discretization of the PDEs.
In the case of geometric multigrid methods, we employ a hierarchy of nested meshes $\{{\pazocal{T}}_{\ell} \}_{\ell=0}^L$, which encapsulate the computational domain $\Omega$.
Through the following, we use the subscript $\ell=0, \ldots, L$ to denote a level, where $L$ denotes the finest level and $0$ denotes the coarsest level.
We denote the number of unknowns on a given level as $\{n_\ell\}_{\ell=0}^{L}$.

The multigrid method relies on three main ingredients.
Firstly, a set of transfer operators is required to pass the information between the subsequent levels of the multilevel hierarchy.
Secondly, suitable smoothers are needed to damp the high-frequency components of the error associated with a given level $\ell$.
Finally, an appropriate coarse level solver is required to eliminate the low-frequency components of the error.
As the JFNK methods are inherently matrix-free, these ingredients have to be adapted, such that they give rise to a matrix-free multigrid preconditioner.

\vspace{0.3cm}
\noindent
\textbf{Transfer Operators:}
In the standard multigrid method, the interpolation~$\boldsymbol{I}_{\ell-1}^{\ell}: \mathbb{R}^{n_{\ell-1}} \to \mathbb{R}^{n_{\ell}}$ and
restriction~$\boldsymbol{R}_{\ell}^{\ell-1}:\mathbb{R}^{n_{\ell}} \to \mathbb{R}^{n_{\ell-1}}$ operators are employed to prolongate the correction to a finer level and
restrict the residual to a coarser level, respectively.
The presented multigrid method requires an evaluation of the action of a Jacobian on a vector on all levels of the multilevel hierarchy.
Therefore, the current Newton iterate also has to be transferred to the coarser levels.
To this aim, we employ a projection operator $\boldsymbol{P}_\ell^{\ell-1}: \mathbb{R}^{n_{\ell}} \to \mathbb{R}^{n_{\ell-1}}$.
In our numerical experiments, we use $\boldsymbol{R}_{\ell}^{\ell-1} := (\boldsymbol{I}_{\ell-1}^\ell)^\top$ and $\boldsymbol{P}_\ell^{\ell-1} = 2^{-d}(\boldsymbol{I}_{\ell-1}^\ell)^\top $, where $d$ denotes the spatial dimension in which the problem is defined.
The scaling factor $2^{-d}$ in the definition of the projection operator $\boldsymbol{P}_\ell^{\ell-1}$ is added to ensure that the constant functions are preserved when projecting them from a fine space to a coarse space.

\vspace{0.3cm}

\noindent
\textbf{Smoothers:}
We utilize the three-level Chebyshev semi-iterative method~\cite{davydov_matrix-free_2020}, as its implementation does not require explicit matrix representation.
This method is convergent if all eigenvalues of the Jacobian lie within a bounded interval.
Our aim here is to reduce only the high-frequency components of the error associated with a given level $\ell$.
Therefore, we focus on the interval $[0.06\lambda_{\ell},1.2\lambda_{\ell}]$, where~$\lambda_{\ell}$ is an estimated largest eigenvalue of the Jacobian on the level~$\ell$.
We estimate the eigenvalue~$\lambda_{\ell}$ at the beginning of each Newton iteration.
More precisely, we employ the Power method, which we terminate within $30$ iterations or when the difference between the subsequent estimates is lower than $10^{-2}$.
As an initial guess for the Power method, a random vector is provided at the first Netwon step.
While for the subsequent Newton steps, we utilize the eigenvector associated with the largest eigenvalue, obtained during the previous eigenvalue estimation process, as an initial guess.

\vspace{0.3cm}

\noindent
\textbf{The coarse level solver:}
In the traditional multigrid method, a direct solver is used to eliminate the remaining low-frequency components of the error on the coarsest level.
In the Jacobian-free framework, we replace the direct solver with a Krylov-subspace method, e.g.,~CG method.
However, to obtain an accurate solution, a large number of iterations may be required.
To reduce the amount of work, we employ a preconditioner based on the limited memory BFGS (L-BFGS) quasi-Newton method~\cite{morales_automatic_2000}.
The L-BFGS preconditioner is created during the very first call to the CG method by storing a few secant pairs.
Following~\cite{morales_automatic_2000}, we collect the secant pairs using the uniform sampling method, which allows us to capture the whole spectrum of the Jacobian.

By design, the CG method is suitable for solving the symmetric positive definite systems.
When solving the non-convex problems, the arising linear systems might be indefinite, which can render the CG method ineffective.
To ensure the usability of the CG method, we propose a few modifications.
Firstly, we terminate the iteration process, as soon as the negative curvature is encountered~\cite{jorgenocedal2000-04-27}.
At this point, we also compute the Rayleigh quotient, given as $\lambda_{c} = \Big(\frac{\boldsymbol{p}^\top \boldsymbol{A} \boldsymbol{p}}{ \boldsymbol{p}^\top \boldsymbol{p}} \Big)$, which gives an estimate of the eigenvalue encountered at the current iterate (that will be also negative).
Secondly, we shift the whole spectrum of the Jacobian by adding a multiple of identity, given as ${\boldsymbol{A}_s = \boldsymbol{A} + (-\lambda_c) \boldsymbol{I}}$, where $\boldsymbol{I}$ denotes an identity matrix.
The shifting strategy is applied recursively, until the modified $\boldsymbol{A}_s$ becomes positive definite.
Please note, the application of the $\boldsymbol{A}_s$ to a vector can be evaluated trivially in the Jacobian-free framework.
The shifting parameter $\gamma$ has to be chosen to be large enough that we do not require many shifting iterations and it has to be small enough that the $\lambda_{\min}(\boldsymbol{A}_s) \approx -\lambda_{\min}(\boldsymbol{A})$.

The multigrid algorithm equipped with the shifting strategy is described in Algorithm~\ref{alg:kothari_h_mini_17_mg}.

\begin{algorithm}[t]
  \caption{Jacobian-free Multigrid - $V(\nu_1,\nu_2)$-cycle}
  \label{alg:kothari_h_mini_17_mg}
  \SetKwData{Left}{left}\SetKwData{This}{this}\SetKwData{Up}{up}
  \SetKwComment{Comment}{$\triangleright$\ }{}
  \SetKwFunction{Union}{Union}\SetKwFunction{FindCompress}{FindCompress}
  \SetKwInOut{Input}{input}
  \SetKwInOut{Output}{output}
  Function: $\boldsymbol{s}_\ell \mapsfrom$ MG($\boldsymbol{x}_\ell^{(k)}, F(\boldsymbol{x}_\ell^{(k)}), \boldsymbol{b}_\ell, \ell$)\\
  $\boldsymbol{s}_{\ell}    \mapsfrom \boldsymbol{0} $                                                                \Comment*[r]{Initialize correction}
  \eIf{$\ell \neq 0$}{
    $\boldsymbol{s}_{\ell}    \mapsfrom $\ Smoother$(\boldsymbol{s}_\ell,\boldsymbol{x}_\ell^{(k)}, F(\boldsymbol{x}_\ell^{(k)}), \boldsymbol{b}_\ell, \nu_1)$  \Comment*[r]{Pre-smoothing}
$\boldsymbol{r}_{\ell-1}  \mapsfrom \boldsymbol{R}_{\ell}^{\ell-1} (\boldsymbol{b}_\ell - J(\boldsymbol{x}_\ell^{(k)})\boldsymbol{s}_\ell)$                                \Comment*[r]{Restrict the residual}
    $\boldsymbol{x}_{\ell-1}^{(k)} \mapsfrom \boldsymbol{P}_{\ell}^{\ell-1} \boldsymbol{x}_\ell^{(k)} $                         \Comment*[r]{Restrict Newton iterate}
    $\boldsymbol{c}_{\ell-1}  \mapsfrom $\ MG$(\boldsymbol{x}_{\ell-1}^{(k)}, F(\boldsymbol{x}_{\ell-1}^{(k)}), \boldsymbol{r}_{\ell-1}, \ell-1)$      \Comment*[r]{Recursion}
    $\boldsymbol{s}_{\ell}    \mapsfrom \boldsymbol{s}_\ell + \boldsymbol{I}_{\ell-1}^{\ell} \boldsymbol{c}_{\ell-1}$                    \Comment*[r]{Update the correction}
    $\boldsymbol{s}_{\ell}    \mapsfrom $\ Smoother$(\boldsymbol{s}_\ell,\boldsymbol{x}_\ell^{(k)}, F(\boldsymbol{x}_\ell^{(k)}), \boldsymbol{b}_\ell, \nu_2)$  \Comment*[r]{Post-smoothing}
  }
  {
    $\lambda_{c+} \mapsfrom 0$ \Comment*[r]{Initialize shifting parameter}
    $\boldsymbol{s}_0,\lambda_{c} \mapsfrom $\ CG$(\boldsymbol{s}_0,\boldsymbol{x}_0^{(k)}, F(\boldsymbol{x}_0^{(k)}), \boldsymbol{r}_0, \lambda_{c+} ,\nu_\ast)$   \Comment*[r]{Coarse level solver}
    \While{$\lambda_{c} < 0$}{
      $\lambda_{c+}  \mapsfrom \gamma \min( \lambda_c, \lambda_{c+} )$ \Comment*[r]{Update shifting parameter}
$\boldsymbol{s}_0,\lambda_{c} \mapsfrom $\ CG$(\boldsymbol{s}_0,\boldsymbol{x}_0^{(k)}, F(\boldsymbol{x}_0^{(k)}), \boldsymbol{r}_0, \lambda_{c+} ,\nu_\ast)$   \Comment*[r]{Shifted CG solver}
    }
  }
\end{algorithm}

\section{Numerical Experiments}
\label{sec:kothari_h_mini_17_numerics}
We investigate the performance of the proposed MG preconditioner through three examples.
We note, for these examples the analytical form of the Jacobian can be computed, but following the JFNK methods, we restrict ourselves from using this information or assembling the Jacobian on the coarsest level.
We use discretize then optimize approach, where the discretization is done with the first order FE method.\\
\noindent\textbf{Bratu:}
Let us consider a domain $\Omega:= (0,1)^2$.
The solution of Bratu problem is obtained by solving the following energy minimization problem:
\begin{equation}
  \begin{aligned}
    \min_{u \in H^1(\Omega)} \Psi_B(u) & = \int_{\Omega} \frac{1}{2} \|\nabla u \|^2 - \lambda \exp(u) \ d\boldsymbol{x}, \\
    \text{such that}  \  \  u          & = 0 \  \text{on} \ \Gamma,
  \end{aligned}
\end{equation}
where we choose $\lambda=5$ and $\Gamma = \partial \Omega$ denotes the boundary.
In our experiments, the mesh ${\pazocal{T}}_{0}$ is triangular and consists of $25$ elements in each direction.

\noindent
\textbf{Minimal Surface:}
We consider again a domain $\Omega:= (0,1)^2$.
This experiment aims to find the surface of minimal area described by the function $u$ by solving the following convex minimization problem:
\begin{equation}
  \begin{aligned}
    \min_{u \in H^1(\Omega)}  \Psi_M(u) & = \int_\Omega \sqrt{(1 + \|\nabla u \|^2)} \ d\boldsymbol{x}, \\
    \text{such that} \   \  u           & = 0 \   \quad \quad \ \ \   \text{on} \ \Gamma_{D_1} ,  \\
    u                                    &= x(1-x)\ \text{on} \ \Gamma_{D_2},
  \end{aligned}
  \label{eq:kothari_h_mini_17_min_surf}
\end{equation}
where, $\Gamma_{D_1}= \{[0,y) \cup [1,y) \}$ and $\Gamma_{D_2}= \{(x,0] \cup (x,1]\}$.
We consider mesh $\pazocal{T}_0$ as in the previous example.

\noindent
\textbf{Hyperelasticity:}
At the end, we investigate a finite strain deformation of a beam, $\Omega = (0,10)\times(0,1)\times(0,1)$, with the rotational deformation applied on the boundaries
$\Gamma_{D_1}=\{0\}\times[0,1]\times[0,1]$, and $\Gamma_{D_2}=\{10\}\times[0,1]\times[0,1]$.
We consider Neo-Hookean material model, and seek for the displacement field $\boldsymbol{u}$ by solving the following non-convex minimization problem:
\begin{equation}
  \begin{aligned}
    \min_{\boldsymbol{u} \in [H^1(\Omega)]^3}  \Psi_N(\boldsymbol{u}) & = \int_{\Omega} \frac{\mu}{2}(I_C-3) - \mu(\ln(J)) + \frac{\lambda}{2} (\ln(J))^2 \ d \boldsymbol{x}, \\
    \text{such that} \  \  \boldsymbol{u}                    & =  \vec{0}  \   \quad  \ \   \text{on} \ \Gamma_{D_1} ,                                       \\
    \boldsymbol{u}                                           & =  \boldsymbol{u}_2  \   \quad    \text{on} \  \ \Gamma_{D_2} ,                                        \\
  \end{aligned}
\end{equation}
where $\boldsymbol{u}_2 = (0, 0.5(0.5+(y-0.5)\cos(\pi/6)-(z-0.5)\sin(\pi/6)-y),0.5(0.5+(y-0.5)\sin(\pi/6)+(z-0.5)\cos(\pi/6)-z))$.
Here, $J:= \det(\boldsymbol{F})$ denotes the determinant of the deformation gradient~$\boldsymbol{F}:= \boldsymbol{I} +\nabla \boldsymbol{u}$.
The first invariant of the right Cauchy-Green tensor is computed as $I_C:= \text{trace}(\boldsymbol{C})$, where $\boldsymbol{C}=\boldsymbol{F}^\top\boldsymbol{F}$.
For our experiment, the Lam\'e parameters $\mu = \frac{E}{2(1+\nu)}$ and $\lambda = \frac{E \nu}{(1+\nu)(1-2\nu)}$ are obtained by setting the value of Young's modulus $E = 10$ and Poisson's ratio $\nu=0.3$.
On the coarse level, the domain is discretized using hexahedral mesh, denoted as ${\pazocal{T}}_0$, with $10$ elements in $x$-directions and $1$ elements in $y$ and $z$ directions.

\vspace{0.3cm}

\noindent
\textbf{Setup for the solution strategy:} We solve the proposed numerical examples using the inexact JFNK (IN) method with a cubic backtracking line-search algorithm~\cite{dennis1983numerical}.
At each IN iteration, the search direction is required to satisfy $\|J(\boldsymbol{x}^{(k)}) \delta \boldsymbol{x}^{(k)} +  F(\boldsymbol{x}^{(k)}) \| \leq \eta^{(k)}\|F(\boldsymbol{x}^{(k)})\|$, where ${\eta^{(k)}=\min(0.5,\|F(\boldsymbol{x}^{(k)})\|)}$.
The algorithm terminates if $\|F(\boldsymbol{x}^{(k)})\|< 10^{-6}$.
We solve $J(\boldsymbol{x}^{(k)}) \delta \boldsymbol{x}^{(k)} = -F(\boldsymbol{x}^{(k)})$, using three different solution strategies: the CG method without any preconditioner (CG), the CG method with L-BFGS preconditioner (CG-QN), and the CG method with the multigrid preconditioner (CG-MG).
The L-BFGS preconditioner is constructed during the first inexact Newton iteration by storing $20$ secant pairs.
The V-cycle MG preconditioner performs $5$ pre-smoothing and $5$ post-smoothing steps.
On the coarse level, we use the CG-QN method with the spectral shift, which is activated only if the negative curvature is encountered.
We employ a shifting parameter $\gamma=5$, in Algorithm~\ref{alg:kothari_h_mini_17_mg}.
The coarse level solver terminates if $\|\boldsymbol{r}_0\| \leq 10^{-12}$, or if the maximum number of iterations, given by the number of unknowns, is reached.

The performance of all solution strategies is evaluated for increasing problem size on successively finer refinement levels.
The refinement levels are denoted by $L0, L1\ldots, L5$, where $L0$ denotes the coarse level, equipped with mesh $\pazocal{T}_0$. The number of levels in the multilevel hierarchy is increased with the refinement level, e.g.,~MG employs $2$ levels for the $L1$ refinement level and $6$ levels for the $L5$ refinement level.
We assess the performance of the methods by measuring the number of required gradient evaluations (GE).
In multilevel settings, the number of effective gradient evaluations is computed as
$\mathrm{GE}= \sum_{\ell=0}^L 2^{-d(L-\ell)} \mathrm{GE}_{\ell}$,
where $\mathrm{GE}_{\ell}$ denotes the number of gradient calls on a given level $\ell$.

We note, the discretization of the minimization problem is performed using the finite element framework libMesh~\cite{kirk2006libmesh}, while the presented solution strategies are implemented as a part of the open-source library UTOPIA~\cite{utopiagit}.

\vspace{0.3cm}

\noindent
\textbf{Influence of different preconditioners on the performance of the JFNK method:}
Table~\ref{tab:kothari_h_mini_17_gradient_evals} and~\ref{tab:kothari_h_mini_17_num_iters} illustrate the performance of the IN method with different linear solvers.
As we can see, for the smaller problems ($L1, L2$), the IN method with the CG and the CG-QN outperforms the IN method with the CG-MG method.
However, as the problem size increases, the IN method with CG-MG is significantly more efficient than with CG or CG-QN.
For instance, for the Bratu example and $L5$ refinement level, the CG-MG method outperforms the other methods by an order of magnitude.

The nonlinearity of the Bratu problem is not affected by the problem size and therefore the number of IN iterations remains constant for all refinement levels.
We can also observe that the behavior of the CG-MG method is level-independent.
The number of required gradient evaluations is therefore bounded after few refinements, as the cost of the coarse level solver becomes negligible.
The same behavior can not be observed for the minimal surface problem, as this problem is strongly nonlinear and the nonlinearity of the problem grows with increasing problem size.
Due to this reason, the number of IN iterations and the total gradient evaluations also increases for the minimal surface problem.
However, we note, that increase is more prevalent for IN method equipped with the CG or the CG-QN methods than with the CG-MG method.

For the hyperelasticity example, the stored energy functional is non-convex hence the negative curvature is quite often encountered on the coarse level.
We notice that with increasing problem size, the negative curvature is encountered fewer times.
As a consequence, a huge amount of coarse level gradient evaluations is required to shift the spectrum of the Jacobian for smaller problems.
Therefore, the average number of gradient evaluations per CG-MG decreases as the problem size increases, as we can observe in Table~\ref{tab:kothari_h_mini_17_num_iters}.
Nevertheless, IN method equipped with the CG-MG outperforms the CG and the CG-QN methods, see Table~\ref{tab:kothari_h_mini_17_gradient_evals}.
Interestingly, the use of the L-BFGS preconditioner is less effective, as in the first IN iteration, the CG method terminates before the whole spectrum of the Jacobian can be captured.
\begin{table*}[t]
  \begin{subtable}[]{1\linewidth}\centering
    \setlength{\tabcolsep}{0.5em}
    \smallskip
    \begin{tabular}{ c " r | r | r " r | r | r " r | r| r}
      \multirow{2}{*}{Levels}              &
      \multicolumn{3}{c"}{Bratu}           &
      \multicolumn{3}{c"}{Minimal surface} &
      \multicolumn{3}{c}{Hyperelasticity}                                                                       \\
      \cline{2-10}
                                           & CG   & CG-QN & CG-MG & CG   & CG-QN & CG-MG & CG   & CG-QN & CG-MG \\\thickhline
      \rowcolor{kothari_Gray}
      $L1$                                 & 176  & 107   & 264   & 360  & 229   & 596   & 467  & 546   & 868   \\
      $L2$                                 & 367  & 233   & 253   & 835  & 501   & 567   & 626  & 655   & 372   \\
      \rowcolor{kothari_Gray}
      $L3$                                 & 767  & 476   & 244   & 2009 & 1170  & 662   & 1349 & 1464  & 426   \\
      $L4$                                 & 1582 & 1097  & 239   & 3544 & 2201  & 782   & 1971 & 1954  & 733   \\
      \rowcolor{kothari_Gray}
      $L5$                                 & 3377 & 2345  & 238   & 6154 & 4316  & 931   & --   & --    & --
    \end{tabular}
  \end{subtable}
\caption{The number of total gradient evaluations required in inexact JFNK method.}
  \label{tab:kothari_h_mini_17_gradient_evals}
\end{table*}
\begin{table*}[t]
  \begin{subtable}[]{1\linewidth}\centering
    \setlength{\tabcolsep}{0.4em}
    \begin{tabular}{ c " r | r | r " r | r | r " r | r| r}
      \multirow{2}{*}{Levels}              &
      \multicolumn{3}{c"}{Bratu}           &
      \multicolumn{3}{c"}{Minimal surface} &
      \multicolumn{3}{c}{Hyperelasticity}                                                                                        \\
      \cline{2-10}
                                           & \# IN & \# CG-MG & \# AGE & \# IN & \# CG-MG & \# AGE & \#IN & \# CG-MG & \# AGE \\\thickhline
      \rowcolor{kothari_Gray}
      $L1$                                 & 3     & 7        & 39.32   & 6     & 13       & 45.85   & 9    & 28       & 34.66   \\
      $L2$                                 & 3     & 9        & 28.25   & 7     & 18       & 31.51   & 5    & 15       & 25.12   \\
      \rowcolor{kothari_Gray}
      $L3$                                 & 3     & 9        & 27.10   & 8     & 25       & 26.50   & 5    & 20       & 20.95   \\
      $L4$                                 & 3     & 9        & 26.61   & 9     & 32       & 24.45   & 5    & 39       & 18.76   \\
      \rowcolor{kothari_Gray}
      $L5$                                 & 3     & 9        & 26.53   & 9     & 41       & 22.79   & --   & --       & --
    \end{tabular}
  \end{subtable}
  \caption{The total number of inexact JFNK iterations ($\#$ IN),  the total number of CG-MG iterations ($\#$ CG-MG), and the average number of gradient evaluations per total linear iteration ($\#$ AGE).}
  \label{tab:kothari_h_mini_17_num_iters}
\end{table*}
\begin{table*}[t]
  \begin{subtable}[]{1\linewidth}\centering
    \setlength{\tabcolsep}{0.55em}
    \begin{tabular}{ c " r | r | r " r | r | r " r | r| r}
      \multirow{2}{*}{Levels}    &
      \multicolumn{3}{c"}{CG}    &
      \multicolumn{3}{c"}{CG-QN} &
      \multicolumn{3}{c} {Shifted CG-QN}                                                                         \\
      \cline{2-10}
                                 & \# IN & \# CG-MG & \#GE  & \# IN & \# CG-MG & \#GE  & \# IN & \# CG-MG & \#GE \\\thickhline
      \rowcolor{kothari_Gray}
      $L1$                       & 9     & 3010     & 51094 & 9     & 130      & 2057  & 9     & 28       & 868  \\
      $L2$                       & 5     & 16       & 44866 & 5     & 1017     & 14814 & 5     & 15       & 372  \\
      \rowcolor{kothari_Gray}
      $L3$                       & 5     & 21       & 1265  & 5     & 26       & 512   & 5     & 20       & 426  \\
      $L4$                       & 6     & 39       & 733   & 6     & 39       & 733   & 5     & 39       & 733
    \end{tabular}
  \end{subtable}
  \caption{
    The total number of inexact JFNK iterations ($\#$ IN),  the total number of CG-MG iterations ($\#$ CG-MG), and the total number of gradient evaluations ($\#$~GE) with
 CG, CG-QN, and shifted  CG-QN methods.
    The experiment was performed for the hyperelasticity example.
  }
  \label{tab:kothari_h_mini_17_neohookean}
\end{table*}

\vspace{0.3cm}
\noindent
\textbf{Effect of the coarse level solver on the performance of the multigrid:}
Due to the non-convexity of the stored energy function, for the hyperelasticity problem, it becomes essential to shift the spectrum of the Jacobian on the coarse level to retain the performance of the multigrid preconditioner.
If only CG or CG-QN method is used, the total number of effective gradient evaluations blows up, as we can see in Table~\ref{tab:kothari_h_mini_17_neohookean}.
This is due to the fact, that the coarse level solver (CG/CG-QN method) terminates as soon as the negative curvature is encountered.
Therefore, the low-frequency components of the error are not eliminated and the multigrid preconditioner becomes unstable.
In contrast, if we employ the shifting strategy, the multigrid preconditioner becomes stable and the total number of the gradient evaluations grows in proportion with the number of required linear iterations.

In conclusion, the performed experiments demonstrate that the proposed Jacobian-free multigrid is a robust and stable preconditioner when applied to problems of various types.
Additionally, we observe level-independence behavior, if the nonlinearity or non-convexity of the problem is not influenced by the discretization parameter.
\vspace{0.1cm}

\footnotesize{
\noindent
\textbf{Acknowledgements:}
The authors would like to thank the Swiss National Science Foundation for their support through the project and the Deutsche Forschungsgemeinschaft (DFG) for their support in the SPP 1962 “ Stress-Based Methods for Variational Inequalities in Solid Mechanics: Finite Element Discretization and Solution by Hierarchical Optimization [186407]”.
Additionally, we would also like to gratefully acknowledge the support of Platform for Advanced Scientific Computing (PASC)  through projects FASTER: Forecasting and Assessing Seismicity and Thermal Evolution in geothermal Reservoirs.
}


\begin{thebibliography}{99.}

\bibitem{bastian_matrix-free_2019}
Bastian, P., M{\"u}ller, E.H., M{\"u}thing, S., Piatkowski, M.: Matrix-free
  multigrid block-preconditioners for higher order discontinuous {Galerkin}
  discretisations.
\newblock Journal of Computational Physics \textbf{394}, 417--439 (2019)

\bibitem{davydov_matrix-free_2020}
Davydov, D., Pelteret, J.P., Arndt, D., Kronbichler, M., Steinmann, P.: A
  matrix-free approach for finite-strain hyperelastic problems using geometric
  multigrid.
\newblock International Journal for Numerical Methods in Engineering
  \textbf{121}(13), 2874--2895 (2020)

\bibitem{dennis1983numerical}
Dennis, J.E., Schnabel, R.B.: Numerical methods for nonlinear equations and
  unconstrained optimization.
\newblock Classics in Applied Math \textbf{16} (1983)

\bibitem{kirk2006libmesh}
Kirk, B.S., Peterson, J.W., Stogner, R.H., Carey, G.F.: libmesh: a c++ library
  for parallel adaptive mesh refinement/coarsening simulations.
\newblock Engineering with Computers \textbf{22}(3-4), 237--254 (2006)

\bibitem{knoll_jacobian-free_2004}
Knoll, D.A., Keyes, D.E.: Jacobian-free {Newton}{\textendash}{Krylov} methods:
  a survey of approaches and applications.
\newblock Journal of Computational Physics \textbf{193}(2), 357--397 (2004)

\bibitem{mavriplis_assessment_2002}
Mavriplis, D.J.: An {Assessment} of {Linear} {Versus} {Nonlinear} {Multigrid}
  {Methods} for {Unstructured} {Mesh} {Solvers}.
\newblock Journal of Computational Physics \textbf{175}(1), 302--325 (2002)

\bibitem{may_scalable_2015}
May, D.A., Brown, J., Le~Pourhiet, L.: A scalable, matrix-free multigrid
  preconditioner for finite element discretizations of heterogeneous {Stokes}
  flow.
\newblock Computer Methods in Applied Mechanics and Engineering \textbf{290},
  496--523 (2015)

\bibitem{morales_automatic_2000}
Morales, J.L., Nocedal, J.: Automatic {Preconditioning} by {Limited} {Memory}
  {Quasi}-{Newton} {Updating}.
\newblock SIAM Journal on Optimization \textbf{10}(4), 1079--1096 (2000)

\bibitem{jorgenocedal2000-04-27}
Nocedal, J., Wright, S.: Numerical Optimization.
\newblock Springer (2000)

\bibitem{utopiagit}
Zulian, P., Kopani{\v c}{\'a}kov{\'a}, A., Nestola, M.C.G., Fink, A., Fadel, N., Rigazzi, A., Magri, V., Schneider, T., Botter, E., Mankau, J., Krause, R.: {U}topia: {A} {C}++ embedded domain specific language for scientific
  computing. {G}it repository.
\newblock https://bitbucket.org/zulianp/utopia (2016)
\end{thebibliography}
\end{document}